\documentclass[11pt]{amsart}

\usepackage{amsthm,amssymb,amsmath}
\usepackage{newtxtext,newtxmath}
\usepackage[
    paperwidth=5.5in,
    paperheight=8.5in,
    top=0.35in,
    left=0.35in,
    right=0.35in,
    bottom=0.5in
]{geometry}
\usepackage[dvipsnames,svgnames,table]{xcolor}
\usepackage[colorlinks=true,linkcolor=RoyalBlue,urlcolor=RoyalBlue,citecolor=PineGreen]{hyperref}
\usepackage{microtype}
\usepackage[nameinlink,capitalise]{cleveref}

\newcommand{\R}{\mathbb{R}}
\newcommand{\D}{\mathcal{D}}
\newtheorem{theorem}{Theorem}
\newtheorem{lemma}[theorem]{Lemma}

\theoremstyle{remark}
\newtheorem{remark}[theorem]{Remark}

\title[A refinement of Szemer\'edi's distinct-distance bound for points in general position]
{A refinement of Szemer\'edi's distinct-distance bound for points in general position}
\author{P\'eter \'Agoston}
\address{P\'eter \'Agoston\\Charles University, Czech Republic}
\email{agoston@kam.mff.cuni.cz}
\author{Krishnendu Bhowmick}
\address{Krishnendu Bhowmick\\TU Chemnitz, Germany}
\email{krish.combi@gmail.com}
\date{}

\makeatletter
\renewcommand{\@settitle}{%
  \begin{center}
    \normalfont\fontsize{14}{16}\selectfont\bfseries
    \MakeUppercase{\@title}\par
  \end{center}}
\renewcommand{\@setauthors}{%
  \begin{center}
    \hypersetup{urlcolor=black}
    \normalfont\fontsize{11}{13}\selectfont
    \begin{tabular}[t]{@{}c@{}}
      {\fontsize{12}{14}\selectfont P\'eter \'Agoston\textsuperscript{*}}\\[2pt]
      Charles University\\
      Czech Republic\\
      \href{mailto:agoston@kam.mff.cuni.cz}{agoston@kam.mff.cuni.cz}
    \end{tabular}\hfill
    \begin{tabular}[t]{@{}c@{}}
      {\fontsize{12}{14}\selectfont Krishnendu Bhowmick}\\[2pt]
      TU Chemnitz\\
      Germany\\
      \href{mailto:krish.combi@gmail.com}{krish.combi@gmail.com}
    \end{tabular}
  \end{center}}
\renewcommand{\@setaddresses}{}
\makeatother
\begin{document}

\begin{abstract}
We prove that every collinear triple free set of $n\ge4$ points in the Euclidean plane  determines at least $\lfloor n/3\rfloor+1$ distinct
distances. This refines the classical bound $\lceil(n-1)/3\rceil$ of
Szemer\'edi that is based on counting isosceles triangles. The refinement is based on treating the diameters separately: they only appear on the convex hull, so for them, we can use the planar diameter theorem and Altman's theorem on
distances in convex polygons.
\end{abstract}
\maketitle
\thispagestyle{empty}
\begingroup
  \renewcommand{\thefootnote}{\fnsymbol{footnote}}
  \footnotetext[1]{Supported by grant no.~23-04949X of the Czech Science Foundation (GA\v{C}R).}
\endgroup

\section{Introduction}

For a finite set $P\subset\R^2$, let
\[
 \D(P)=\{|p-q|:p,q\in P,\ p\ne q\},
 \qquad d(P)=|\D(P)|.
\]
The study of the minimum possible value of $d(P)$ goes back to
Erd\H os~\cite{Erdos46}. In his original paper, Erdős proved that if by $d(n)$, we denote the minimum value of $d(P)$ over all point sets in the plane, then $d(n)=O\left(n/\sqrt{log{n}}\right)$, which can be obtained using an integer grid. This bound is conjectured to be tight. The celebrated work of Guth and Katz \cite{GK} using the Elekes--Sharir framework \cite{ES} almost resolved this problem by proving that any point set $P \subset \mathbb R^2$ determines $\Omega\left ( n/ \log n \right )$ distinct distances. 

When no three points of $P$ are collinear,
Szemer\'edi, whose argument is recorded by
Erd\H os~\cite[pp.~100--101]{Erdos75}, proved
\begin{equation}\label{eq:classical}
 d(P)\ge\left\lceil\frac{|P|-1}{3}\right\rceil.
\end{equation}
The vertices of a regular $n$-gon determine $\lfloor n/2\rfloor$
distances, and it is conjectured that \(d(P)\geq \lfloor |P|/2\rfloor\) for any collinear triple free set \(P\subset\R^2\); see also~\cite[Problem~6]{Sheffer}.

The purpose of this note is to improve the integer rounding in
\eqref{eq:classical}.

\begin{theorem}\label{thm:main}
Let $P$ be a set of $n\ge4$ points in $\R^2$, no three of which are
collinear. Then
\[
 d(P)\ge\left\lfloor\frac n3\right\rfloor+1
       =\left\lceil\frac{n+1}{3}\right\rceil.
\]
\end{theorem}

Compared with \eqref{eq:classical}, the bound increases by one for
$n\equiv0,1\pmod3$ and is unchanged for $n\equiv2\pmod3$. The weaker
rounding $d(P)\ge\lceil n/3\rceil$ is recorded in the online
note~\cite[Section~1]{Working1082}, where the equality case in
Szemer\'edi's argument is excluded using the diameter graph. Thus the
additional case addressed here is $n=3d(P)$.

The proof uses two classical results: a planar $s$-point set determines
at most $s$ pairs at its diameter, and the vertices of a strictly convex
$h$-gon determine at least $\lfloor h/2\rfloor$ distances. If $d(P)\le n/3$,
the isosceles-triangle count and the diameter graph force at least
$4n/5$ points to be vertices of the convex hull. Altman's theorem then
reduces the argument to a six-point configuration, which we exclude
directly.

\section{Proof of the main theorem}

Fix a set $P$ of $n\ge4$ points with no three collinear. Write
$d=d(P)$, and list its distances as
\[
 \ell_1<\ell_2<\cdots<\ell_d=\Delta.
\]
For $p\in P$ and $1\le j\le d$, put
\[
 r_{p,j}=\bigl|\{q\in P\setminus\{p\}:|p-q|=\ell_j\}\bigr|.
\]
The classes are indexed by the global distance set, so some $r_{p,j}$
may be zero. At every point,
\begin{equation}\label{eq:row-sum}
 \sum_{j=1}^d r_{p,j}=n-1.
\end{equation}
\begin{lemma}\label{lem:count}
With this notation,
\begin{equation}\label{eq:second-moment}
 \sum_{p\in P}\sum_{j=1}^d r_{p,j}^2\le3n(n-1).
\end{equation}
\end{lemma}

\begin{proof}
Let
\[
 T=\sum_{p\in P}\sum_{j=1}^d\binom{r_{p,j}}2.
\]
This counts isosceles triangles with a distinguished apex; in
particular, an equilateral triangle is counted three times. For a fixed
unordered base pair, every possible apex belongs to its perpendicular
bisector. Since no three points of $P$ are collinear, there are at most
two such apices. Hence
\[
 T\le2\binom n2=n(n-1).
\]
Using \eqref{eq:row-sum}, we obtain
\[
 \sum_{p,j}r_{p,j}^2=2T+\sum_{p,j}r_{p,j}
 \le2n(n-1)+n(n-1),
\]
which proves the lemma.
\end{proof}

For comparison, Cauchy--Schwarz gives
\[
 \frac{n(n-1)^2}{d}
 \le\sum_{p,j}r_{p,j}^2\le3n(n-1),
\]
recovering $d\ge(n-1)/3$. To improve the rounding, we retain the
integrality of the distance-class sizes and treat the diameter
separately.

We recall the geometric facts used below. The planar diameter theorem
states that a set of $s$ points determines at most $s$ unordered pairs
at its diameter; see Erd\H os~\cite[Theorem~3]{Erdos46}. Altman proved
that a strictly convex $h$-gon determines at least $\lfloor h/2\rfloor$
distances~\cite{Altman63}; see also~\cite[Theorem~3]{Altman72}.

Let $H$ be the set of vertices of the convex hull of $P$, and put
$h=|H|$. Since no three points are collinear, $H$ is the vertex set of
a strictly convex polygon. Altman's theorem implies
\begin{equation}\label{eq:altman}
 \left\lfloor\frac h2\right\rfloor\le d,
 \qquad\text{and hence}\qquad h\le2d+1.
\end{equation}
Every endpoint of a diameter of $P$ belongs to $H$. Indeed, if
$|p-q|=\Delta$, then $P$ is contained in the closed disk with centre
$q$ and radius $\Delta$. The tangent to that disk at $p$ is a
supporting line that meets the disk only at $p$, so $p$ is a
convex-hull vertex.

For $p\in P$, write $a_p=r_{p,d}$ for its degree in the diameter graph.
The preceding observation and the diameter theorem applied to $H$ give
\begin{equation}\label{eq:diameter-degrees}
 a_p=0\quad(p\notin H),
 \qquad \sum_{p\in H}a_p\le2h.
\end{equation}
Here the diameter of $H$ is $\Delta$, because both endpoints of every
diameter pair of $P$ lie in $H$.

\begin{lemma}\label{lem:hull}
If $d\le n/3$, then $h\ge4n/5$.
\end{lemma}

\begin{proof}
For every integer $r$, the inequality $(r-3)(r-4)\ge0$ gives
\[
 r^2\ge7r-12.
\]
Apply this to the $r_{p,j}$'s for $j<d$. At a
point $p\in P$, their sizes sum to $n-1-a_p$, so
\[
 \sum_{j=1}^d r_{p,j}^2
 \ge a_p^2+7(n-1-a_p)-12(d-1).
\]
This also holds when $d=1$, with the sum over $j<d$ empty. Summing
over $p$ gives
\begin{equation}\label{eq:refined-moment}
 \sum_{p,j}r_{p,j}^2
 \ge n(7n+5-12d)+\sum_{p\in P}(a_p^2-7a_p).
\end{equation}
For every integer $a$, one also has $(a-2)(a-3)\ge0$, or
\[
 a^2-7a\ge-2a-6.
\]
The points of $P\setminus H$ have $a_p=0$. Applying the last inequality only
on $H$, and then using \eqref{eq:diameter-degrees}, we obtain
\begin{align*}
 \sum_{p\in P}(a_p^2-7a_p)
 &=\sum_{p\in H}(a_p^2-7a_p)\\
 &\ge-2\sum_{p\in H}a_p-6h\ge-10h.
\end{align*}
Combining this with \eqref{eq:refined-moment} and \cref{lem:count}
gives
\[
 3n(n-1)\ge n(7n+5-12d)-10h,
\]
and therefore
\[
 10h\ge4n(n-3d)+8n\ge8n,
\]
where the last inequality uses $d\le n/3$.
\end{proof}

\begin{proof}[Proof of \cref{thm:main}]
Suppose that $d\le n/3$. By \cref{lem:hull} and \eqref{eq:altman},
\[
 \frac{4n}{5}\le h\le2d+1\le\frac{2n}{3}+1.
\]
Consequently $n\le15/2$, and hence $4\le n\le7$. The integer bounds
\[
 \left\lceil\frac{4n}{5}\right\rceil\le h
 \le2\left\lfloor\frac n3\right\rfloor+1
\]
exclude $n=4,5,7$: their respective lower bounds for $h$ are $4,4,6$,
while the upper bounds are $3,3,5$. Thus $n=6$. The classical bound
\eqref{eq:classical} and the assumption $d\le n/3$ give $d=2$.

In this case \eqref{eq:altman} gives $h\le5$, so there is a point
$c\in P\setminus H$. Such a point is not incident with a diameter.
Writing $R=\ell_1$, all five points of $P\setminus\{c\}$ therefore
lie on the circle $C$ with centre $c$ and radius $R$.

Fix $p\in P\setminus\{c\}$. The circles centred at $p$ with radii $R$
and $\Delta$ each meet $C$ in at most two points, since their centres
differ from $c$. The four points of $P\setminus\{c,p\}$ all lie on
these two circles. Hence exactly two of them are at distance $R$ from
$p$. On $C$, the points at distance $R$ from $p$ are precisely its
rotations about $c$ through angles $\pi/3$ and $-\pi/3$. Both
rotations therefore belong to $P\setminus\{c\}$.

It follows that this five-point set is invariant under rotation
through $\pi/3$ about $c$. Every orbit of that rotation on $C$ has
six elements, a contradiction. Thus $d>n/3$, proving the theorem.
\end{proof}

\begin{remark}
The restriction $n\ge4$ is necessary: an equilateral triangle
determines one distance. The bound is attained for $n=4,5,6,7$ by the
regular $n$-gon. It leaves unchanged the leading coefficient $1/3$ in
the classical lower bound.
\end{remark}

\section{Concluding remarks}

Since obtaining the result in this paper, the result has been substantially improved, see \cite{Bhowmick}.

\section*{Acknowledgements and AI disclaimer}
The authors thank Dömötör Pálvölgyi for helpful suggestions. The authors also thank the organizers of KAMAK 2026, where the authors met.

    This work was carried out in collaboration with ChatGPT. The authors take full responsibility for the contents of the paper.

\bibliography{references}{}
\bibliographystyle{plain}
\end{document}